\documentclass[10pt,leqno]{article}
\AtBeginDvi{} 
\usepackage{amssymb}
\usepackage{amscd}
\usepackage{amsmath}
\usepackage{amsthm}
\usepackage{mathrsfs}
\usepackage{graphics}
\usepackage{graphicx}
\def\overset#1#2{{\mathrel{\mathop {{#2}_{}}\limits^{#1}}}}
\def\underset#1#2{{\mathrel{\mathop {{}_{} {#2}}\limits_{{#1}_{}}}}}
\def\upplim_#1{\underset{#1}{\overline\lim}\;}
\def\lowlim_#1{\underset{#1}{\underline\lim}\;}

\newtheorem{claim}[equation]{\indent \rm Claim}%

\newtheorem{defn}[equation]{\indent{\it Definition}\rm }
\newtheorem{exam}[equation]{\indent \rm {\it Example}}
\newtheorem{llist}[equation]{\rm\hskip2pt }
\newtheorem{lem}[equation]{Lemma}
\newtheorem{prop}[equation]{Proposition}
\newtheorem{rmk}[equation]{\indent \rm {\it Remark}}

\newtheorem{thm}[equation]{Theorem}

\newcommand{\alg}{\mathrm{alg}}
\newcommand{\C}{{\mathbf{C}}}

\newcommand{\codim}{{\mathrm{codim}}}

\newcommand{\iso}{\cong}

\newcommand{\lto}{\longrightarrow}

\newcommand{\ponec}{{\mathbf{P}^1(\mathbf{C})}}

\newcommand{\Q}{{\mathbf{Q}}}
\newcommand{\R}{{\mathbf{R}}}

\newcommand{\st}{{\mathrm{St}^0}}

\newcommand{\tor}{\mathrm{tor}}

\newcommand{\Z}{\mathbf{Z}}
\newcommand{\zar}{\mathrm{Zar}}
\newenvironment{pf}{\par{\it Proof.\,}}{\qed\vskip+2pt}
\numberwithin{equation}{section}

\title{
An application of the value distribution theory for semi-abelian varieties
to problems of Ax-Lindemann and Manin-Mumford types
}
\date{The University of Tokyo}
\author{
Junjiro Noguchi\thanks{Research supported in part by Grant-in-Aid
 for Scientific Research (C) 15K04917. \hfill\break
\hbox{\quad} MSC2010: 11J95, 32H30, 03C64  \hfill\break
\hbox{\quad} Key words: Ax-Lindemann, Manin-Mumford, big Picard,
Nevanlinna theory, semi-abelian variety \hfill\break
}
}
\begin{document}
\AtBeginDvi{} 
%
\setlength{\baselineskip}{13pt}
\parskip+2.5pt
\maketitle
\begin{abstract}
The aim of this paper is to prove a theorem of Ax-Lindemann type
for complex semi-abelian varieties as an application of a big Picard
theorem proved by the author in 1981, and then apply it to prove a theorem
of classical Manin-Mumford Conjecture for semi-abelian varieties,
which was proved by M. Raynaud 1983, M.~Hindry 1988,
$\ldots$, and Pila-Zannier 2008 by a different method from others,
which is most relevant to our's.
The present result might be a first instance of a {\em direct connection}
between the value distribution theory of holomorphic maps and
the arithmetic (Diophantine) theory over algebraic number fields,
 while there have been many {\em analogies} between them.
\end{abstract}

\section{Introduction and  main results}\label{int}

The purpose of this paper is to prove a theorem of Ax-Lindemann type
for complex semi-abelian varieties as an application of a big Picard
theorem obtained in \cite{n81}, and then apply it to prove a theorem
of Manin-Mumford type\footnote{\, It was called the Manin-Mumford
Conjecture
and proved by M. Raynaud \cite{ra83}, M. Hindry \cite{hi88},
$\ldots$, Pila-Zannier \cite{pz08}; cf.\  \cite{pz08}, Introduction,
S. Lang \cite{la91}, Chap.\ I~\S6,
and e.g., P. Tzermias \cite{tz00} for surveys of the
 Manin-Mumford Conjecture.}
for the distribution of torsion points on a subvariety of a semi-abelian variety
defined over a number field, combined with extending a part of the arguments
in Pila-Zannier \cite{pz08} for abelian varieties (cf.\ \S\ref{prmm}).

\begin{thm}[Ax-Lindemann type]\label{ax-lind}
Let $\exp : \C^n \to A$ be an exponential map of a complex semi-abelian
variety $A$. Let $V \subset \C^n$ be a complex irreducible affine
 algebraic subvariety with the restricted map $\exp|_V: V \to A$.
Then  the Zariski closure $X(\exp|_V)$
of the image of $\exp|_V$ in $A$ is a translate of a complex
 semi-abelian subvariety of $A$.
\end{thm}

\begin{thm}[Manin-Mumford type]\label{mm}
Let $X \subset A$ be a proper algebraic reduced subvariety of
a semi-abelian variety $A$ defined over an algebraic number field.
Then, the Zariski closure $\overline{X \cap A_\tor}^\zar$ of torsion points
on $X$ is a finite union of translates of 
semi-abelian subvarieties by torsion points on $X$.
\end{thm}

In view of the above two theorems, it is naturally led
to study a denseness property of the
value-distribution of $\exp|_V$ in Theorem \ref{ax-lind}.
In fact, for an ample divisor $D$
on a projective compactification of $A$, we will prove that
$D \cap \exp (V)$ is Zariski dense in $D \cap X(\exp|_V)$
(cf.\ Theorem \ref{ddense}).

In the course of the proofs we will introduce a new notion of
``{\em strictly transcendental}'' for holomorphic maps into semi-abelian
varieties (see Definition \ref{stt}). By making use of a Big Picard
Theorem due to \cite{n81} we prove the image structure stated in Theorem
\ref{ax-lind} for strictly transcendental holomorphic maps into
semi-abelian varieties (see Theorem \ref{stim}).
We then prove that the map $\exp|_V$ in Theorem \ref{ax-lind} is
 strictly transcendental (Proposition~\ref{strans}).

The proof of Theorem \ref{mm} roughly consists of two parts: The first
is a decomposition of the torsion points on $X$ to
 an algebraic part and its compliment, done by the arguments due to Pila-Zannier
\cite{pz08}, Theorem 2.1, Step 1,
 being extended to the semi-abelian case.
The second part is the application of Theorem  \ref{ax-lind} to the algebraic
part of torsion points, and for its compliment we use
the lower and upper estimates
due to Masser and Pila-Wilkie (cf.\ \cite{pz08}, \S3) to deduce
the statement of Theorem \ref{mm}.

 In the last \S\ref{ex-st} we give some example of a transcendental but
 not strictly transcendental map into an abelian variety, and study
its image structure.

The present result might be a first instance of a {\em direct connection}
between the value distribution theory of holomorphic maps and
the arithmetic (Diophantine) theory over algebraic number fields,
although there have been many {\em analogies} between them.
It is the point of interest of this paper to observe that the direct
connection is provided by the theory of o-minimal structures
in model theory (Pila-Wilkie \cite{pw06}; this approach evolved from
Bombieri-Pila \cite{bp89} (cf.\  Zannier \cite{za12}).
 This aspect should be of some interest for
 both sides of the theories.

{\it Acknowledgment.} The present paper has grown up through
the interesting and helpful discussions with Professors Pietro Corvaja
and Umberto Zannier since Conference
 ``Specialization Problems in Diophantine Geometry''%
, Cetraro (Italy), 9--14 July 2017.
In particular, the proof of Theorem \ref{mm} in \S\ref{prmm} is based
on their arguments suggested through the discussions.
It is a great pleasure of the author to express his deep gratitude to
both of them.

\section{Big Picard and Ax-Lindemann}

{\bf (i) Big Picard. }
 Let $\Delta^m$ be the unit polydisk of $\C^m$ with center at $0$ and let
$E \subset \Delta^m$ be a thin complex analytic subset. Let $Y$ be a
reduced complex space with a compactification $\bar Y$, a compact
complex space with reduced structure.
 Let $f: \Delta^m \setminus E \to Y$  be a holomorphic map.
If there is a meromorphic map $\bar f : \Delta^m \to \bar Y$
with the restriction to $\Delta^m \setminus E$,
 $\bar f |_{(\Delta^m \setminus E)}=f$,
  we say that $E$ is {\it removable} for $f$ or $f$ is
(meromorphically) extendable over $E$: Otherwise, $f$ is said to be
{\it transcendental}.

We denote by $X(f)$ the Zariski closure of the image of $f$ in $Y$
with respect to the Zariski topology on $\bar Y$.
Note that $X(f)$ is irreducible.

\begin{rmk}[Hartogs extension theorem]\label{hart}\rm
Let $Y$ be quasi-projective algebraic.
If $\codim\, E \geq 2$, $E$ is always removable for any holomorphic map
 $f: \Delta^m \setminus E \to Y$.
 Therefore, as the removability is
 concerned with quasi-projective $Y$,
 it suffices to deal with smooth $E$ of codimension $1$.
\end{rmk}

Let $A$ be a complex semi-abelian variety with a smooth projective
compactification $\bar A$, and let $X \subset A$ be a non-empty
 complex algebraic subset. We denote by $\st(X)$ the connected component
 of $0$ of the group $\{a \in A: a+X=X\}$, and simply call it
the {\it stabilizer} of $X$. Then, $\st(X)$ is an algebraic subgroup of
$A$ and a semi-abelian subvariety by itself.
Note that $X$ is of general type if and only if $\st(X)=\{0\}$.

We recall:
\begin{thm}[\cite{n81}, Corollary (4.7)]\label{bpic}
Assume that $X$ is an irreducible algebraic subvariety of $A$.
Let $f: \Delta^m \setminus E \to X \hookrightarrow A$ be a holomorphic
 map.
If $\st(X(f))=\{0\}$, then $f$ is extendable over $E$; in the other
 words, if $f$ is transcendental, then $\dim \st(X(f))>0$.
\end{thm}

\begin{defn}[strictly transcendental]\label{stt}\rm
A transcendental holomorphic map $f:\Delta^m \setminus E \to A$ is said
 to be {\it strictly transcendental} if for every semi-abelian subvariety
$B \subset A$ the composite $q_B \circ f: \Delta^m \setminus E \to A/B$
with the quotient map $q_B: A \to A/B$ is either constant or
transcendental.
\end{defn}

\begin{thm}\label{stim}
Let $f: \Delta^m \setminus E \to A$ be a
strictly transcendental holomorphic map.
Then, $X(f)=f(c)+\st(X(f))$ with $c \in \Delta^m \setminus E$;
i.e., $X(f)$ is a translate of a semi-abelian subvariety of $A$.
\end{thm}

\begin{pf}
Set $B=\st(X(f))$ and $g=q_B \circ f: \Delta^m \setminus E \to A/B$.
It suffices to show that $g$ is constant.
Otherwise, $g$ would be transcendental by the assumption for $f$.
Since $\st(X(g))=\{0\}$, it follows from Theorem \ref{bpic} that
$g$ is extendable over $E$; this is a contradiction.
\end{pf}

\smallskip\noindent
{\bf (ii) Ax-Lindemann. }
We keep the notation used above. 
Here we always assumes that ``varieties'' are irreducible with reduced structure.
Let
\[
 \exp: \C^n \to A
\]
be an exponential map.
A map from a complex algebraic variety into another complex
algebraic variety is called a {\it transcendental} map if it is
{\em not} a rational map.
We then define a
 {\it strictly transcendental map} from a complex algebraic variety
into $A$ as in Definition \ref{stt}.

\begin{prop}\label{strans}
Let $V \subset \C^n$ be a positive dimensional irreducible
complex algebraic subvariety of $\C^n$.
 Then, the restricted map
$\exp |_V: V \to A$ of $\exp$ to V is strictly transcendental.
\end{prop}
\begin{pf}
Let $B \subsetneq A$ be a complex semi-abelian subvariety.
Suppose that
$$q_B \circ \exp |_V: V \lto A_1=A/B$$
 is a non-constant rational map. We shall deduce
a contradiction.

There is  an algebraic curve  $C \subset V$ (irreducible
 and $\dim C=1$) such that the restriction
$$f:= q_B \circ \exp|_C: C \lto A_1$$
 is non-constant  rational.

(a) Case of $B=\{0\}$:
 There is an exact sequence
\begin{equation}\label{ex-sa}
0 \to (\C^*)^t \to A \to A_0 \to 0,
\end{equation}
where $A_0$ is an abelian variety. Set $L=\exp^{-1} (\C^*)^t$.
Then $L$ is a $t$-dimensional vector subspace of $\C^n$.
Let $p_0: \C^n \to \C^n/L \iso \C^m ~ (m=n-t)$ and $q_0: A \to A_0$
be the quotient maps. Then $\exp$ naturally induces an
exponential map
\[
 \exp_0: \C^m \to A_0.
\]
We are going to infer 
\begin{claim}\label{clm}
 $q_0 \circ \exp\, (C)$ is a point.
\end{claim}

Suppose that it is not the case.
Then, the Zariski closure $C_0$ of the image
$p_0(C)$ in $\C^m$ is an  algebraic curve in $\C^m$.
Let $\omega_0$ be a flat K\"ahler metric on $A_0$ such that
$\exp_0^*\omega_0=\alpha$, where
 $\alpha=\sum_{j=1}^m \frac{i}{2\pi} dz_j \wedge d\bar{z}_j$
with the natural coordinate system $(z_1, \ldots, z_m)$ of $\C^m$.
Since $\exp|_C$ is rational, $\exp_0|_{C_0}$ is non-constant rational.
 Then,  the Zariski closure  $W:=\overline{\exp_0 (C_0)}^\zar$
in $A_0$ is an algebraic curve in $A_0$ with
\begin{equation}
\int_{\exp_0 (C_0)} \omega_0 =\int_{W} \omega_0
=M < \infty.
\end{equation}
If $k$ denotes the degree of the rational map
$\exp_0|_{C_0}: C_0 \to W$, we have
\begin{equation}\label{finite}
\int_{C_0} \alpha  = kM.
\end{equation}
Let $B(r) \subset \C^m$ be an open ball of radius $r (>0)$ with center at
a point $a \in C_0$ and set $C_0(r)=C_0 \cap B(r)$.
Wirtinger's inequality implies
\[
 \int_{C_0(r)} \alpha \geq \nu(a; C_0) r^2,
\]
where $\nu(a; C_0) (\geq 1)$ denotes the order of $C_0$ at $a$.
Letting $ r \to \infty$, we have a contradiction to \eqref{finite}.
Therefore, Claim~\ref{clm} follows.

Thus, we have a non-constant rational map after a translation:
\[
\exp|_C:  C (\subset L \iso \C^t) \to (\C^*)^t,
\]
which is a restriction of exponential map
\[
 (\zeta_1, \ldots, \zeta_t) \in \C^t \to \left (e^{\zeta_1}, \ldots,
e^{\zeta_t} \right) \in (\C^*)^t.
\]
Let $\bar{C}$ be the closure of $C$ in $(\ponec)^t \supset (\C^*)^t$.
Here we write $\ponec=\C \cup\{\infty\}$.
Then there is a point $b=(b_1, \ldots, b_t) \in \bar{C} \setminus C$
with some $b_j=\infty$. Since $e^{\zeta_j}$ has an isolated \textit{essential}
singularity at $\zeta_j=\infty$, $\exp |_C$ cannot be rational.

(b) Case of general $B$: Set
\begin{align*}
 F &=\exp^{-1} B,\\
 p' &: \C^n \to  \C^n/F \iso \C^{n'}, \\
 q' &: A \to A/B = A'.
\end{align*}
Let $\exp': \C^{n'} \to A'$ be the naturally induced exponential from
 $\exp: \C^n\to A$.
Let $C'$ be the Zariski closure of $p'(C)$ in $\C^{n'}$.
Then, it follows from the assumption that $\pi'|_{C'}: C' \to A'$
would be a non-constant rational map: This contradicts
to what was proved in (a) above.
\end{pf}

{\em Proof of Theorem \ref{ax-lind}}:  This is now immediate by
Proposition \ref{strans} and Theorem \ref{stim}. \qed

\begin{exam}\rm Here we give a simple example of a strictly
 transcendental map.
Set
\begin{align*}
 \exp &: (z,w) \in \C^2 \to (e^z, e^w) \in (\C^*)^2 =A ,\\
 V &= \{(z, w)\in \C^2: zw=1\}=\{(z, 1/z): z \in \C^* \} \iso \C^*.
\end{align*}
Then $V$ is affine algebraic, 
\[
\exp_V : (z, 1/z) \in V \lto (e^z, e^{1/z}) \in A
\]
is strictly transcendental, and $X(\exp|_V)=A$.
\end{exam}

\section{Proof of Theorem \ref{mm}}\label{prmm}
Let $K$ be a number field over which $A$ is defined.
Let $$\exp: \C^n \lto A$$ be an exponential map.
The reduction of $\exp^{-1} X \cap A_\tor$ to the algebraic part
is done in parallel to the proof of Pila-Zannier \cite{pz08}, \S3,
and then we apply Theorem \ref{ax-lind} to conclude the theorem.

The proof is done by induction on $\nu=\dim X \geq 0$.
If $\nu=0$, it is trivial. Suppose that the case of
$\dim X \leq \nu-1 (\nu \geq 1)$ holds. Let $\dim X=\nu$.

We consider the stabilizer $\st(X)$.

(a) The case of $\dim \st(X) >0$. We set the quotient map
$q: A \to A_1=A/\st(X)$ and $X_1=q(X)=X/\st(X)$.
Since $\dim X_1 <\nu$, the induction hypothesis implies that
there are at most finitely many semi-abelian subvarieties
$B_j \subset A_1, 1 \leq j \leq l$, and torsion points $b_j \in X_1 \cap A_{1\tor}$
such that
\[
 X_1 \cap A_{1\tor}=\bigcup_{j=1}^l (b_j+B_{j \tor}).
\]
Taking any elements $a_j \in q^{-1} b_j \cap X \cap A_\tor, 1 \leq j \leq l$,
we have
\[
X \cap A_{\tor} \subset \bigcup_{j=1}^l (a_j + (q^{-1}B_j)_\tor).
\]
Since $X \supset a_j+ q^{-1} B_j$, the opposite inclusion above also holds; thus,
\[
\overline{X \cap A_{\tor}}^\zar = \bigcup_{j=1}^l (a_j + q^{-1}B_j).
\]

(b) The case of $\st(X)=\{0\}$. Note that $X$ is of (log) general type.
We put $\Lambda=\exp^{-1}\{0\} ~(\subset \C^n)$ and $Z=\exp^{-1} X$
which is a $\Lambda$-periodic analytic subset of $\C^n$.
We consider a basis of $\Lambda$ as $\Z$-module, which is also
a basis of the real vector space $\R \Lambda$ of dimension
$d=2n -t$. Thus we have identifications
\[
 \Z^d \iso \Z \Lambda \subset \R \Lambda \iso \R^d.
\]
The torsion points on $X$ lift to points on
 $\Q \Lambda \iso \Q^d \subset \R^d$.
We consider the set $\check Z= Z \cap \Q^d$ of the preimage
of all torsion points on $X$  and its restriction
$\check Z_1:=\check Z \cap [0,1]^d$ of $\check Z$
to the closed fundamental domain $[0, 1]^d$.

For a subset $W \subset \Q^d$ and a real number $T \geq 1$, we denote
by $N_W(T)$ the number of rational points in $W$ whose denominator is
at most $T$. 

By Pila-Wilkie \cite{pw06} there is a so-called algebraic part
$\check Z_1^{\alg}$ of $\check Z_1$,
 satisfying the following properties:
\begin{enumerate}
\item
We have
\begin{equation}\label{ap}
\check Z_1^{\alg} = \bigcup_{V \subset Z} (V \cap \{t=(t_j) \in
 [0,1)^d: t_j \in \Q\}),
\end{equation}
where $V$ runs over all positive dimensional
affine algebraic subsets of $\C^n$,
contained in $Z$ (cf.\ \cite{pz08}, Proposition 2.1\footnote{\,It
 is noted that the periodicity condition of Proposition 2.1
 is not used in the proof, so that it can be applied to our $Z$.}).
\item
For every $\varepsilon >0$ there is a positive constant
$c_1 (=c_1(\check Z_1, \varepsilon)$,
  depending on $\check Z_1$ and $\varepsilon$)
 such that
\begin{equation}\label{uest}
N_{\check Z_1'}(T)
 \leq c_1 T^\varepsilon, \quad T \geq 1,
\end{equation}
where $\check Z_1':= \check Z_1 \setminus \check Z_1^{\alg}$,
called the transcendental part of $\check Z_1$ (\cite{pw06}).
\end{enumerate}

We first analyze the algebraic part $\check Z_1^{\alg}$.
If $\check Z_1^{\alg}=\emptyset$, it is done; we assume
$\check Z_1^{\alg} \not = \emptyset$.
Let $V$ be as in \eqref{ap}. It follows from Theorem \ref{ax-lind}
that $X(\exp|_V) ~(\subset X)$ is a translate of a positive dimensional
 semi-abelian subvariety of $A$.
Let $Y$ denote the union of all translates $a+B'$ ($a \in X$) of
positive dimensional semi-abelian subvarieties $B' \subset A$ such that
$a+B' \subset X$.
It follows from \cite{n81}, Lemma (4.1) (cf.\ also \cite{nw14}, \S5.6.4,
and \cite{pz08}, p.\,160)
that $Y\,(\subset X)$ is an algebraic subset
of dimension $<\nu$ and for every irreducible component $Y'$ of $Y$,
$\dim \st(Y')>0$.
We may assume that $Y$ is defined over $K$.
Applying the induction hypothesis to $Y$, we have
finitely many positive dimensional semi-abelian subvarieties
 $B'_j \subset A, 1 \leq j \leq l'$, and $P_j \in Y \cap A_\tor$
such that
\[
 Y \cap A_\tor=\bigcup_{j=1}^{l'} (P_j+B'_{j \tor}).
\]
Therefore, we have
\begin{align*}
 \exp\,(Z_1^\alg) = \bigcup_{j=1}^{l'} (P_j+B'_{j \tor}), \quad
\overline{\exp\,(\check Z_1^{\alg})}^\zar =
\bigcup_{j=1}^{l'} (P_j+B'_j).
\end{align*}

We may assume that all $P_j$ and $B'_j$ are defined over $K$
after a finite extension of $K$. Then, we have:
\begin{llist}[Invariance]\label{coninv}\rm
 For a torsion point
$P \in \exp\,(\check Z_1^\alg)$ (resp.\ $\exp\,(\check Z_1')$),
 all of its conjugates over $K$
lie in $\exp\,(\check Z_1^\alg)$ (resp.\ $\exp\,(\check Z_1')$).
\end{llist}

To analyze the part $\check Z_1'$ we prepare:
\begin{lem}\label{lest}
Let $P \in A$ be a torsion point of exact order $N$.
Then, there is a number $\rho >0$ depending only on $\dim A$
such that
\begin{equation}\label{lest2}
 [K(P): K] \geq c_2 N^\rho, \quad T \geq 1,
\end{equation}
where $c_2=c_2(A, K)$ is a positive constant with dependence as signed.
\end{lem}

\begin{pf}
For the semi-abelian variety $A$
we have the following exact sequence
\[
 0 \to \mathbf{G}_m^t \to A \, \overset{\pi}{\lto} \, A_0 \to 0,
\]
where $\mathbf{G}_m^t$ is an algebraic torus and
$A_0$ is an abelian variety;
 we may assume that all of the above morphisms and
algebraic groups are defined over $K$.
Then, $P_0:=\pi(P)$ is a torsion point of $A_0$, whose exact order is
denoted by $N_0$. Note that $N_0 P \in \mathbf{G}_m^t$ and
$N=N_1 N_0$, where $N_1$ is the order of $N_0P ~(\in \mathbf{G}_m^t)$.
Then,
\[
 c_3 \varphi(N_1) \leq [ K(N_0 P): K] \leq  [K(P): K],
\]
where $c_3=c_3(t, K)$ is a positive constant and $\varphi$ is the
Euler function. It is known that
\[
 \varphi (N_1) \sim \frac{N_1}{\log \log N_1}.
\]
If $N_1 \geq \sqrt{N}$, the proof is finished.
Otherwise, we have $N_0 \geq \sqrt{N}$. We then apply Masser's
estimate (\cite{m84}) for the following second inequality:
\[
 c_3 N^{\rho'/2} \leq c_3 N_0^{\rho'} \leq [K(P_0): K] \leq   [K(P): K] ,
\]
where $c_3=c_3(A_0, K)$ is a positive constant.
Thus, we have the lower estimate \eqref{lest2}.
\end{pf}

\smallskip
{\it Finish of the proof of Theorem \ref{mm}}\,:
Let now $P\in \exp\,(\check Z_1')$ be a torsion point of exact order
$N$.
It follows from Invariance~\ref{coninv} that
$N_{\check Z_1'}(N) \geq [K(P): K]$.
By Lemma \ref{lest}
we see that
\begin{equation}\label{lest3}
 N_{\check Z_1'} (N) \geq c_2 N^\rho.
\end{equation}
On the other hand, we have by \eqref{uest}
\begin{equation}\label{uest2}
 N_{\check Z_1'} (N) \leq  c_1 N^\varepsilon.
\end{equation}
Taking $\varepsilon = \rho/2$ with $\rho$ in \eqref{lest3},
we conclude the boundedness of $N$ by \eqref{lest3} and \eqref{uest2}.
Therefore, $\check Z_1'$ is a finite set.
\qed

\section{Distribution of $\exp (V)$ on divisors}\label{vdis}

Let $A$ be a semi-abelian variety with exponential map
\[
 \exp: \C^n \lto A,
\]
and let $V \subset \C^g$ be an irreducible affine algebraic subvariety.
Let $\bar A$ be a projective compactification of $A$.
Because of the results of the previous sections, it might be of some
interest to look at the actual value-distribution of $\exp|_V (V)$
in its Zariski closure $X(\exp|_V)$.

We first deal with a transcendental holomorphic map
$f: \Delta^* \to A$ from a punctured disk $\Delta^*$ into $A$.
In \cite{nwy08} we dealt with entire holomorphic maps from the whole
plane $\C$ into $A$. Combining the method of \cite{n81}
with the result and the arguments explored in
 \cite{nwy08} and \cite{nw14} we see that the results
of \cite{nwy08} (and \cite{nw14}, Chap.\ 6) hold
 for transcendental holomorphic maps from
$\Delta^*$ into $A$. In particular, we have (cf., also \cite{cn12},
Theorem 5.2)
\begin{thm}\label{dense}
Let $f: \Delta^* \to A$ be a transcendental holomorphic map and let
$D$ be an effective ample divisor on $\bar A$ such that
 $D \not\supset f(\Delta^*)$.
Then $f(\Delta^*) \cap X(f) \cap D$ is Zariski dense in
$X(f) \cap D$ (where $D$ stands also for the support of $D$).
\end{thm}

By making use of this we prove:
\begin{thm}\label{ddense}
Let $V \subset \C^n$ be an irreducible complex affine algebraic subvariety,
and let $D$ be an effective ample divisor on $\bar A$.
Then, the intersection  $D \cap \exp(V)$ is Zariski dense in
$ D \cap X(\exp|_V)$.
\end{thm}
\begin{pf}
Let $\zeta_0 \in V$ be fixed.
We consider a pencil of affine algebraic curves $C_\gamma \subset V,
\gamma \in \Gamma$, passing through $\zeta_0$, such that
$\bigcup_{\gamma} C_\gamma$ contains a non-empty open subset
of $V$ in the sense of differential topology.
By Theorem \ref{ax-lind}  $X(\exp|V)$ and $X(\exp|_{C_\gamma})$
are all translates of semi-abelian subvarieties of $A$ passing
through $\exp\,(\zeta_0)$. Since there are at most countably many
such semi-abelian subvarieties,
 one finds a curve $C_0=C_\gamma$ such that
 $X(\exp|_{C_0})=X(\exp|_V)$.
Then it suffices to show the theorem for $C_0$.
Let $C_1 \to C_0$ be the normalization and let $\bar C_1$ be its
smooth compactification.
Then there is an analytic neighborhood $U ~(\subset \bar C_1)$ of a point $Q$ of
$\bar C_1 \setminus  C_1$ such that $U \setminus \{Q\}$
is biholomorphic to a punctured disk $\Delta^*$.
Then our assertion is immediate by Theorem \ref{dense}.
\end{pf}

\section{An example of a transcendental but not strictly transcendental
 map}\label{ex-st}
Let $\bar{C}$ be a smooth complex  projective algebraic curve of genus $g \geq 1$ and
let $q: \bar{C} \to J(\bar{C})$ be the Jacobian embedding; here, when
$g=1$, we simply take $q: \bar C \to \bar C (=J(\bar C))$ as the identity map.
We set $A_1=J(\bar{C})$, which is an abelian variety of dimension $g$.
Let $Q \in \bar{C}$ be a point and set $C=\bar{C}\setminus \{Q\}$.
Then $C$ is affine algebraic and there is a
finite map $p: C \to \C$.

Let $\exp : \C^g \to A_1$ be an exponential map.
We take a linear embedding $\lambda: \C \to \C^g$
that is in sufficiently generic direction with respect
to the period lattice of $\exp : \C^g \to A_1$.
Then, $X(\exp \circ \lambda )=A_1$.
We put
\begin{equation}\label{tmap}
f: x \in C \lto (q(x), \exp \circ \lambda \circ p (x)) \in 
A_1 \times A_1=: A.
\end{equation}
\begin{prop}\label{non-st}
Let $f: C \to A$ be as above.
\begin{enumerate}
\item
The holomorphic map $f$ is transcendental but not strictly
     transcendental.
\item
If $g \geq 2$, the Zariski closure $X(f)$ of the image $f(C)$ is
{\em not} a translate of an abelian subvariety of $A$.
\item
If $g=1$, $X(f)=A$.
\end{enumerate}
\end{prop}
\begin{pf}
(i) The first half is clear. For the latter, note
 that $ X(f) \subset \bar{C} \times A_1$.
With a subgroup $\{0\}\times A_1 \subset A_1 \times A_1=A$
we consider the quotient map $\mu:A \to A/\{0\}\times A_1 \iso A_1$.
Then, $\mu \circ f=q|_C : C \to  A_1$ is rational.
Therefore, $f$ is not strictly transcendental.

(ii)  Since $\mu (X(f))=q(\bar C)$, $X(f)$ is not
a translate of an abelian subvariety.

(iii) Since $\dim X(f)=1$, or $2$, it suffices to deduce a
 contradiction with assuming $\dim X(f)=1$. If so, $X(f)$ is a translate
of an abelian subvariety of $A$. We consider an effective divisor
$D=\bar C \times \{w\}$ with $w \in \bar C$.
We infer from the definition of $f$ that $X(f) \cap D$ is infinite.
Therefore, $X(f)=D$; this is a contradiction.
\end{pf}

\begin{rmk}\label{paun-sibony}\rm\hskip-3pt\footnote{~
The author is glad to acknowledge that the present work
was inspired by the talks given by J. Pila and E. Ullmo
at the Cetraro Conference 2017 referred at the end of \S\ref{int},
 and the talk of Ullmo drew his attention to the preprint \cite{ps17}.}
P\u{aun}-Sibony \cite{ps17} deals with a similar application
of the Bloch-Ochiai Theorem to the abelian Ax-Lindemann statement
(\cite{ps17}, Theorems 5.2 and 5.6).
But with regard to Proposition \ref{non-st} above,
in Theorems 5.2 and 5.6 of \cite{ps17}
one might be able to have only the
non-triviality of the stabilizer of the Zariski closure
of the image as in Theorem \ref{bpic}, obtained in \cite{n81}
(Corollary (4.7)).
\end{rmk}

\bigskip
\setlength{\baselineskip}{12pt}
\rightline{Graduate School of Mathematical Sciences}
\rightline{University of Tokyo (Emeritus)}
\rightline{Komaba, Meguro-ku, Tokyo 153-8914}
\rightline{Japan}
\rightline{e-mail: noguchi@ms.u-tokyo.ac.jp}
\end{document}